\RequirePackage{fix-cm}
\documentclass[smallextended]{svjour3}       \smartqed  \usepackage{graphicx}
\usepackage{xcolor}
\usepackage{enumerate}
\usepackage{tikz}
\usepackage{pgfplots}
\usepackage{tudacolors}
\usepackage{siunitx}
\usepackage{verbatim}
\usepackage{hyperref}
\usepackage{acronym}
\usepackage{amsmath}
\usepackage{amsfonts}
\usepackage{amssymb}

\renewcommand{\div}{\operatorname{div}}
\newcommand{\eps}{\varepsilon}
\newcommand{\RR}{\mathbb{R}}

\begin{document}

\title{On torque computation in electric machine simulation by harmonic mortar methods}

\titlerunning{Torque computation by harmonic mortar methods}        

\author{Herbert Egger         \and
        Mané Harutyunyan    \and
        Richard L\"oscher     \and
        Melina Merkel\footnotemark[1]      \and
        Sebastian Sch\"ops
}

\institute{Herbert Egger \at
           Johannes Kepler University Linz and Johann Radon Institute for Computational and Applied Mathematics, Altenberger Stra{\ss}e 69, A-4040 Linz\\
           \email{herbert.egger@jku.at,herbert.egger@ricam.oeaw.ac.at}           \and Mané Harutyunyan, Melina Merkel, Sebastian Sch\"ops \at Computational Electromagnetics, TU Darmstadt, Schlossgartenstr. 8, D-64289 Darmstadt\\
  \email{mane.harutyunyan@tu-darmstadt.de, melina.merkel@tu-darmstadt.de, sebastian.schoeps@tu-darmstadt.de}
  \and
  Richard L\"oscher \at 
  Department of Mathematics, TU~Darmstadt, Dolivostra{\ss}e 15, 64293 Darmstadt, Germany\\
  \email{richard.loescher@tu-darmstadt.de}
  \and
  \footnotemark[1]Correspondence: \at 
  \email{melina.merkel@tu-darmstadt.de}
}

\date{Received: date / Accepted: date}

\maketitle

\begin{abstract}
The use of trigonometric polynomials as Lagrange multipliers in the harmonic mortar method enables an efficient and elegant treatment of relative motion in the stator-rotor coupling of electric machine simulation. Explicit formulas for the torque computation are derived by energetic considerations, and their realization by harmonic mortar finite element and isogemetric analysis discretizations is discussed. Numerical tests are presented to illustrate the theoretical results and demonstrate the potential of harmonic mortar methods for the evaluation of torque ripples. 
\end{abstract}

\section{Introduction}

A particular challenge for electric machine simulation is the relative motion of stator and rotor and the computation of quantities of interest depending on the rotation angle, e.g., the torque as a measure for the magneto-mechanic energy conversion.
Various approaches have been proposed to tackle the coupling of subdomain problems on moving geometries, e.g., moving band,  mortar and Lagrange multiplier or discontinuous Galerkin methods \cite{Alotto01,Buffa01,Lange10,Tsukerman95}, all with their specific advantages and shortcomings. 
In addition, different strategies have been discussed for the numerical computation of the torque, e.g., via virtual displacements or the Maxwell stress tensor; see \cite{HenrotteHameyer04} for an overview.

In this paper, we consider harmonic mortar finite element and isogeometric analysis methods proposed in \cite{Bontinck18,DeGersem04}, which are based on finite element or isogeometric analysis approximations of the magnetic problems in the stator and rotor subdomains, coupled by a Lagrange multiplier technique using trigonometric functions.
These methods have several interesting properties: 
\begin{enumerate}[(a)]
\item They are based on the Galerkin approximation of a certain weak formulation of the problem using Lagrange multipliers. This will allow us to give an unambiguous definition of the torque resulting in an exact energy conservation principle also on the discrete level.
\item The computed torque, as a function of the rotation angle, can be shown to be smooth, i.e.,  infinitely differentiable, and therefore no spurious torque ripples are introduced by non-smoothness of the numerical approximation.
\item The use of trigonometric functions for the Lagrange multipliers addition-ally allows for an efficient assembling of the coupling matrices for multiple rotation angles, avoiding the repeated integration of interface terms.
\end{enumerate}
The resulting methods therefore are particularly well-suited for multiquery applications, e.g., the construction of performance maps, the setup of reduced order models, coupled frequency or time domain simulations. This will be illustrated by computation of cogging torque in the numerical tests. 

The remainder of the manuscript is organized as follows: Section \ref{sec:problem} introduces the magnetostatic setting for electric machines in two dimensions. The strong and weak formulations of the problem are presented and the energy balance is used to derive the torque. The discretization of the problem is given in Section \ref{sec:galerkin} using the Galerkin approach. Section \ref{sec:harmonicmortar} discusses the harmonic mortar method and its implementation. Numerical results for the torque computation using the harmonic mortar method are demonstrated in Section \ref{sec:results}. Finally the paper is concluded with a discussion in Section \ref{sec:conclusions}. 

\section{Problem statement}\label{sec:problem}

\begin{figure}
		\centering
\includegraphics{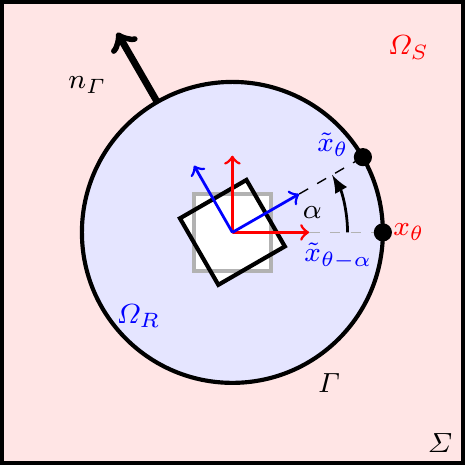}
		\caption{Sketch of the geometry. $x_{\theta}$ denotes a point in the fixed coordinate system (red) of the stator domain $\Omega_S$, whereas $\tilde{x}_{\theta}$ denotes a point in the rotated coordinate system (blue) of the rotor domian $\Omega_R$. The are related by $x_\theta=\tilde{x}_{\theta-\alpha}$ where $\alpha$ is the angle of rotation. }
		\label{fig:geometry}
	\end{figure}

For ease of presentation, we consider a linear magnetostatic setting in two space dimensions. 
The geometric setup represents the cross section of a cylindrical device and consists of a stator domain $\Omega_S$ completely enclosing the rotor domain $\Omega_R$ with common interface $\Gamma = \partial\Omega_S \cap \partial\Omega_R = \partial\Omega_R$; see Figure~\ref{fig:geometry}.
Motivated by the application context, we assume that $\Gamma$ is a circle placed within the air gap and centered around the origin.
Using a standard vector potential formulation \cite{Alotto01,DeGersem04}, the magnetic fields in the two subdomains are described by 
\begin{alignat}{5}
    -\div(\nu_S \nabla a_S) &= j_e \qquad && \text{in } \Omega_S \label{eq:sys1}\\
    -\div(\nu_R \nabla a_R) &= \div M^\perp \qquad && \text{in } \Omega_R \label{eq:sys2}
\end{alignat}
where $a_S$, $a_R$ are the $z$-component of the magnetic vector potential, $j_e$ is the impressed electric source density, and 
$M^{\perp} = (m_{y},-m_{x})$ the rotated magnetization vector.
The equations are formulated in two different coordinate systems attached to the respective subdomains, such that the position of material continuities stays fixed under relative motion.  
We further assume that 
\begin{alignat}{5}
a_S &= 0 && \qquad \text{on } \Sigma \label{eq:sys3}
\end{alignat}
on the outer boundary $ \Sigma = \partial\Omega_S \setminus \Gamma$.
The continuity of the magnetic field and the magnetic vector potential across the interface are encoded in
\begin{alignat}{5}
a_S &= a_R \circ \rho_{-\alpha} \qquad && \text{on } \Gamma, \label{eq:sys4}\\
n_\Gamma \cdot (\nu_S \nabla a_S)  &= n_\Gamma \cdot (\nu_R \nabla a_R) \circ \rho_{-\alpha} \qquad && \text{on } \Gamma, \label{eq:sys5}
\end{alignat}
where $\rho_{\alpha}$ describes the rotation of a point by angle $\alpha$ around the origin.
Furthermore, $n_\Gamma$ denotes the unit normal vector on $\Gamma$ pointing from $\Omega_R$ to $\Omega_S$.
The second coupling condition allows us to introduce 
\begin{alignat}{5}
\lambda :=  n_\Gamma \cdot (\nu_S \nabla a_S) \label{eq:lambda}
\end{alignat}
which amounts to the tangential trace of the magnetic field at the interface. 

\medskip 

\textbf{Dependence on $\alpha$.}
By \eqref{eq:sys4}--\eqref{eq:sys5} the solutions $(a_S,a_R,\lambda)$ of \eqref{eq:sys1}--\eqref{eq:lambda} implicitly depend on the rotation angle $\alpha$, even if the sources $j_e$ and $M$ are independent of $\alpha$, which we assume in the sequel. 
We write $u(\cdot|\alpha)$ to highlight the dependence of a function $u$ on the angle $\alpha$ and denote by 
\begin{align}
u'(\cdot|\alpha) := \tfrac{d}{d\alpha} u(\cdot|\alpha) := \lim_{\eps \to 0} \tfrac{1}{\eps}(u(\cdot|\alpha+\eps) - u(\cdot|\alpha))
\end{align}
the derivative of such a parameter-dependent function $u$ with respect to the parameter $\alpha$.
By differentiation of \eqref{eq:sys4} with respect to $\alpha$, we can see that
\begin{alignat}{5}
a'_S(x_\theta|\alpha) = a'_R(x_{\theta-\alpha}|\alpha) +  \tfrac{d}{d\theta} a_R(x_{\theta-\alpha}|\alpha) \qquad x_\theta \in \Gamma,
\end{alignat}
where $x_{\theta-\alpha} = \rho_{-\alpha}(x_\theta)$ is the point in the rotated coordinate system of $\partial\Omega_R$ corresponding to a point $x_\theta=(r_\Gamma \cos \theta,r_\Gamma \sin \theta)$ on the fixed interface $\Gamma$; see Figure~\ref{fig:geometry}.
Furthermore, $\frac{d}{d\theta} u(x_\theta) = r_\Gamma n_\Gamma^\perp \cdot \nabla u(x_\theta)$ amounts to the tangential derivative of $u$ along $\Gamma$ and $r_\Gamma$ is the radius of $\Gamma$. 

\medskip 

\textbf{Weak formulation.}
Using standard arguments, sufficiently regular solutions of \eqref{eq:sys1}--\eqref{eq:sys5} can be seen to satisfy the variational identities
\begin{alignat}{5}
(\nu_S \nabla a_S, \nabla v_S)_{\Omega_S} + \langle \lambda, v_S\rangle_\Gamma &= (j_e, v_S)_{\Omega_S} \label{eq:var1}\\
(\nu_R \nabla a_R, \nabla v_R)_{\Omega_R} - \langle \lambda, v_R \circ \rho_{-\alpha}\rangle_\Gamma &= (-M^\perp, \nabla v_R)_{\Omega_R} \label{eq:var2}\\
\langle a_S - a_R \circ \rho_{-\alpha}, \mu\rangle_\Gamma &= 0  \label{eq:var3}
\end{alignat}
for all smooth test functions $v_S$, $v_R$, and $\mu$ defined on $\Omega_S$, $\Omega_R$, and $\Gamma$, respectively, with $v_S=0$ on $\Sigma$. 
Using the considerations outlined above, the coupling condition \eqref{eq:var3} can also be rephrased as 
\begin{alignat}{5} 
\langle a'_S - a'_R \circ \rho_{-\alpha}, \mu\rangle_\Gamma &= \langle \tfrac{d}{d\theta} a_R \circ \rho_{-\alpha}, \mu\rangle_\Gamma \label{eq:var3a}
\end{alignat}
which will be used in the following considerations.

\medskip 

\textbf{Energy balance and torque computation.}
As a next step, we introduce the magnetic energy of the system. 
For ease of notation, we assume a unit length $L$  of the device in axial direction. Then 
\begin{align} \label{eq:energy}
    E(a_S,a_R) = \int_{\Omega_S} \tfrac{\nu_S}{2} |\nabla a_S|^2 dx + \int_{\Omega_R} \tfrac{\nu_R}{2} |\nabla a_R|^2 + M^\perp \cdot \nabla a_R \, dx.
\end{align}
Let us note that, via the fields $a_S$, $a_R$, the energy also implicitly depends on the rotation angle $\alpha$. By elementary calculations, we may then compute 
\begin{align*}
\tfrac{d}{d\alpha} E(a_S,a_R) 
&= (\nu_S \nabla a_S,  \nabla a'_S)_{\Omega_S} + (\nu_R \nabla a_R, \nabla a'_R )_{\Omega_R} + (M^\perp, \nabla a'_R)_{\Omega_R} \\
&= (j_e, a'_S)_{\Omega_S} - \langle \lambda,  a'_S\rangle_\Gamma
+ \langle \lambda, a'_R \circ \rho_{-\alpha}\rangle_\Gamma,
\end{align*}
where we used the variational identities \eqref{eq:var1}--\eqref{eq:var2} with $v_S=a'_S$ and $v_R = a'_R$ in the second step. Using the differential form~\eqref{eq:var3a} of the coupling condition, we arrive at the energy balance
\begin{align*}
\tfrac{d}{d\alpha} E(a_S,a_R)  
&= (j_e,a'_S)_{\Omega_S}- \langle \lambda, \tfrac{d}{d\theta} a_R \circ \rho_{-\alpha} \rangle_\Gamma,
\end{align*}
where the first term denotes the electric work required to maintain the electric current $j_e$ in the stator coils (neglecting Ohmic losses), and the second term is the mechanic work required for an infinitesimal rotation of the system. 
These \emph{energy-based} considerations immediately give rise to the following definition of the torque 
\begin{align} \label{eq:torque}
T(\alpha) := \langle \lambda, \frac{d}{d\theta} a_R \circ \rho_{-\alpha} \rangle_\Gamma.
\end{align}
Note that $a_R$ and $\lambda$ and hence $T$ are functions of the rotation angle~$\alpha$. 
In general, the energy and torque have to be scaled by the axial length~$L$. 

\begin{remark}
The above formula can be derived in various ways, see \cite{Bossavit90,Coulomb83,HenrotteHameyer04}. Since the stator and rotor were both assumed to be rigid bodies, the torque can be computed here directly without resorting to the Maxwell stress tensor.
Using rotational symmetry as well as integration by parts along the interface $\Gamma$, we can express the torque alternatively as
\begin{align}
T(\alpha) = -\langle \tfrac{d}{d\theta} \lambda \circ \rho_{\alpha}, a_R\rangle_\Gamma, 
\end{align}
which may be a more convenient representation depending on the particular setting. \end{remark}

\section{Galerkin approximation}\label{sec:galerkin}

For discretization of the problem \eqref{eq:sys1}--\eqref{eq:sys5}, we consider Galerkin approximations of the variational equations \eqref{eq:var1}--\eqref{eq:var3} using appropriate finite-dimensional sub-spaces $V_{h,S} \subset H^1_{\Sigma}(\Omega_S)$, $V_{h,R} \subset H^1(\Omega_R)$, and $M_N \subset H^{-1/2}(\Gamma)$. 
Let us note that the boundary conditions \eqref{eq:sys3} are incorporated explicitly in the definition of $H^1_\Sigma(\Omega_S)=\{v \in H^1(\Omega_S) : v=0 \text{ on } \Sigma\}$ here.  
The discrete problem for a fixed angle $\alpha$ then reads as follows.
\begin{problem} \label{prob:galerkin}
Find $a_{h,S} \in V_{h,S}$, $a_{h,R} \in V_{h,R}$ and $\lambda_N \in M_N$ such that 
\begin{alignat}{5}
(\nu_S \nabla a_{h,S}, \nabla v_{h,S})_{\Omega_S} + \langle \lambda_N, v_{h,S}\rangle_\Gamma &= (j_e, v_{h,S})_{\Omega_S} \label{eq:var1h}\\
(\nu_R \nabla a_{h,R}, \nabla v_{h,R})_{\Omega_R} - \langle \lambda_N, v_{h,R} \circ \rho_{-\alpha}\rangle_\Gamma &= -(M^\perp, \nabla v_{h,R})_{\Omega_R} \label{eq:var2h}\\
\langle \alpha_{h,S} - \alpha_{h,R} \circ \rho_{-\alpha}, \mu_N\rangle_\Gamma &= 0  \label{eq:var3h}
\end{alignat}
for all discrete test functions $v_{h,S} \in V_{h,S}$, $v_{h,R} \in V_{h,R}$, and $\mu_N \in M_N$.
\end{problem}

\textbf{Well-posedness.}
For the following considerations, we assume that $V_{h,S}$, $V_{h,R}$, and $M_N$ are finite-dimensional and that $\nu_S$, $\nu_R$ are uniformly positive and bounded from above. Moreover,the source currents $j_e$ and magnetization vector $M$, as well as the domains $\Omega_S$, $\Omega_R$ are assumed to be sufficiently regular. We can then establish the well-posedness of the discretization scheme as follows.

\begin{lemma} \label{lem:well-posed}
Assume that 
\begin{align} \label{eq:infsup}
\sup_{v_{h,S} \in V_{h,S}} \frac{\langle v_{h,S}, \mu_N\rangle_\Gamma}{\|v_{h,S}\|_{H^1(\Omega_S)}} \ge \beta \|\mu_N\|_{H^{-1/2}(\Gamma)} \qquad \text{for all } \mu_N \in M_N.
\end{align}
Then the discrete variational problem \eqref{eq:var1h}--\eqref{eq:var3h} is uniquely solvable and 
\begin{align*}
    &\|a_S - a_{h,S}\|_{H^1(\Omega_S)} + \|a_R - a_{h,R}\|_{H^1(\Omega_R)} + \|\lambda - \lambda_N\|_{H^{-1/2}(\Gamma)} \\
    &\le C(\beta) \left( \|a_S - v_{h,S}\|_{H^1(\Omega_S)} + \|a_R - v_{h,R}\|_{H^1(\Omega_R)} + \|\lambda - \mu_N\|_{H^{-1/2}(\Gamma)}\right)
\end{align*}
for all $v_{h,S} \in V_{h,S}$, $v_{h,R} \in V_{h,R}$ and $\mu_N \in M_N$.
\end{lemma}
The assertion is a direct consequence of Brezzi's saddle-point theory \cite{Brezzi74}; also see \cite{Buffa01,Egger_2021aa} for application in the current context.
Under the stability condition \eqref{eq:infsup}, the Galerkin approximations obtained by Problem~\ref{prob:galerkin} thus always lead to quasi-optimal error estimates with respect to the chosen approximation spaces. Two particular examples of Galerkin approximations will be discussed below.

\medskip

\textbf{Discrete torque computation and energy balance.}
Mimicking the formula \eqref{eq:torque} for the torque obtained on the continuous level, it seems natural to define the discretized approximation of the torque as
\begin{align} \label{eq:torqueh}
    T_{h,N}(\alpha) := 
    \langle \lambda_h, \tfrac{d}{d\theta} a_{h,R} \circ \rho_{-\alpha}\rangle_\Gamma.
\end{align}
The torque computation for a given angle $\alpha$ thus only requires the solution of one magnetostatic interface problem \eqref{eq:var1h}--\eqref{eq:var3h} for this specific angle.
With the same reasoning as on the continuous level, we obtain the following result.

\begin{lemma} \label{lem:energy_discrete}
Let $(a_{h,S},a_{h,R},\lambda_N)$ denote a solution of \eqref{eq:var1h}--\eqref{eq:var3h}. Then 
\begin{align}
    \tfrac{d}{d\alpha} E(\alpha_{h,S},\alpha_{h,R}) = (j_e,a'_{h,S})_{\Omega_S} - T_{h,N}
\end{align}
with discrete torque $T_{h,N}$ defined by the formula \eqref{eq:torqueh}. 
\end{lemma}
This assertion illustrates that the electro-magneto-mechanic energy balance holds exactly also on the discrete level. 

\medskip 

\textbf{Alternative representation and smooth dependence on $\alpha$.}
If the space $M_N$ of discrete Lagrange multipliers consists of sufficiently smooth functions, we may express the discrete torque equivalently as 
\begin{align} \label{eq:torqueh2}
    T_{h,N} = 
-\langle \tfrac{d}{d\theta} \lambda_N \circ \rho_{\alpha}, a_{h,R}\rangle_\Gamma.
\end{align}
By a recursive argument, one can then show the following result. 
\begin{lemma} \label{lem:smooth}
Let $M_N$ consist of $k-1$ times continuously differentiable functions. 
Then the solution $(a_S,a_R,\lambda_N)$
is $k$ times differentiable with respect to $\alpha$. 
\end{lemma}
A proof of this assertion follows from the algebraic structure of the discretized problem; see the end of the next section. 

\section{Harmonic mortar discretizations.} \label{sec:harmonicmortar}

As a particular choice of a Galerkin approximation, we now consider the following choice of spaces: $V_{h,S}$ and $V_{h,R}$ are standard $H^1$-conforming finite element spaces \cite{Braess} over appropriate triangulations of the stator and rotor domain $\Omega_S$ and $\Omega_R$. Note that no conformity of the mesh across the interface is required. For approximation of the Lagrange multipliers, we consider the space
\begin{align} \label{eq:M} 
M_N = \left\{ \mu_N(x_\theta) = \frac{c_0}{2} + \sum\nolimits_{n=1}^N c_n \cos(n \theta) + d_n \sin(n \theta) : c_n, d_n \in\RR\right\}
\end{align}
of trigonometric polynomials of degree $\le N$. 

\begin{remark} 
As shown in \cite{Egger_2021aa}, the inf-sup condition \eqref{eq:infsup} holds true, if $N$ is chosen appropriately, in particular not too large. By the results of the previous section, the harmonic mortar finite element method then is stable and yields quasi-optimal error estimates. The same conclusions hold for a method based on isogeometric analysis discretization in the stator and rotor domains \cite{Bontinck18,Cottrell_2009aa}.
Following the considerations of the previous section, both methods allow for a torque computation consistent with a discrete energy conservation principle. 
\end{remark}

\textbf{Implementation of the harmonic mortar finite element method.}
In the following, we discuss some details concerning the implementation of the harmonic mortar finite element method outlined above. 
As we will see, the use of trigonometric polynomials for the discrete Lagrange multipliers brings advantages also for the numerical realization.

After choosing appropriate basis functions for the finite element subspaces $V_{h,S}$, $V_{h,N}$ and the trigonometric polynomials in $M_N$, 
the discrete variational problem \eqref{eq:var1h}--\eqref{eq:var3h} can be turned into a linear system
\begin{alignat}{5}
    K_S a_S(\alpha) & &    &+B_S^\top \lambda(\alpha) &&= j_e \label{eq:lin1}\\
    K_R a_R(\alpha) & &    &-B_R(\alpha)^\top \lambda(\alpha) \ &&= j_M \label{eq:lin2}\\
    B_S a_S(\alpha) &-& B_R(\alpha) a_R(\alpha)& &&= 0. \label{eq:lin3}
\end{alignat}
The source vectors $j_e$ and $j_M$ and the matrices $K_S$, $K_R$ and $B_S$ are independent of the angle $\alpha$, whereas $B_R$ and hence also the solution vectors $a_S$, $a_R$ and $\lambda$ depend on the rotation angle. 
For the harmonic mortar, the entries of the coupling matrix $B_R(\alpha)$ can be computed by the following expressions
\begin{align*}
[B_R(\alpha)]_{2n,j} &= \int_0^{2\pi} \cos(n \theta) \, \phi_j(x_{\theta-\alpha}) r_\Gamma d\theta, \qquad n \ge 0,  \\
[B_R(\alpha)]_{2n-1,j} &=  \int_0^{2\pi} \sin(n \theta) \, \phi_j(x_{\theta-\alpha}) r_\Gamma d\theta, \qquad n \ge 1,
\end{align*}
where $\phi_j$ are the basis functions for the Galerkin approximation on $\Omega_R$.

\begin{remark}
By trigonometric summation formulas, we get
\begin{align*}
    [B_R(\alpha)]_{2n,j} &= \cos(n\alpha)[B_R(0)]_{2n,j} -\sin(n\alpha)[B_R(0)]_{2n-1,j}  \\
    [B_R(\alpha)]_{2n-1,j} &= \sin(n\alpha)[B_R(0)]_{2n,j} + \cos(n\alpha)[B_R(0)]_{2n-1,j},
\end{align*}
with smooth coefficients depending on $\alpha$, hence the integrals in the coupling matrix $B_R(\alpha)$ in fact only have to be evaluated for a single angle $\alpha=0$. 
In matrix notation, we may then write 
\begin{align} \label{eq:BRa}
    B_R(\alpha) = R(\alpha) B_R(0)
\end{align}
where $R(\alpha)$ is a block diagonal matrix with $2\times 2$ blocks.
Further, note that the integrals in the definition of $B_R(0)$ consist of products of trigonometric functions with polynomials, which can be computed analytically \cite{Neuman_1981aa}. 
\end{remark}

\textbf{Algebraic energy balance and torque computation.}
On the algebraic level, the energy is given by $E=\frac{1}{2} \|a_S\|^2_{K_S} + \frac{1}{2} \|a_R\|_{K_R}^2 - j_M^\top a_R$ where $\|x\|_M^2 = x^\top M x$ is the square of the norm associated with a symmetric positive semi-definite matrix $M$. Then 
\begin{align*}
    \tfrac{d}{d\alpha} E 
    &= (K_S a_S, a_S') + (K_R a_R, a_R') - (j_M,a_R') \\ 
    &= -(B_S^\top \lambda, a_S') + (j_e, a_S') + (B_R^\top \lambda, a_R') \\
    &= (j_e, a_S') - (B_S a_S' - B_R a_R', \lambda) 
    =(j_e, a_S') - (B_R' a_R, \lambda).
\end{align*}
Hence the discrete torque can be simply computed by the algebraic formula
\begin{align*}
    T_{h,N}(\alpha) = \lambda^\top(\alpha) B_R'(\alpha) a_R(\alpha) =\langle \lambda_h, \tfrac{d}{d\theta}a_{h,R}\circ\rho_{-\alpha}\rangle .
\end{align*}
Note that $B_R'(\alpha)$ can be expressed as $B_R'(\alpha)=R'(\alpha)B_R(0)$ by using \eqref{eq:BRa}. Moreover, the matrix $R'(\alpha)$ is block diagonal with $2 \times 2$ blocks whose entries can be computed analytically; see above.

\medskip

\textbf{Derivatives with respect to $\alpha$ and proof of Lemma~\ref{lem:smooth}.}
By formal differentiation of the above equations, we can see that
\begin{alignat}{5}
    K_S a'_S(\alpha) - B_S^\top \lambda'(\alpha) &= 0 \label{eq:lin1d}\\
    K_R a'_R(\alpha) - B_R(\alpha)^\top \lambda'(\alpha) &= B_R'(\alpha)^\top \lambda(\alpha) \label{eq:lin2d}\\
    B_S a'_S(\alpha) - B_R(\alpha) a'_R(\alpha) &=  B_R'(\alpha) a_R(\alpha). \label{eq:lin3d}
\end{alignat}
Together with \eqref{eq:lin1d}--\eqref{eq:lin3d}, one can see that the solution of \eqref{eq:lin1}--\eqref{eq:lin3} is differentiable with respect to  the angle $\alpha$. 
By recursion, we may further conclude that the torque $T_{h,N}$ of the harmonic mortar method is a smooth function of the angle $\alpha$. The same argument provides a proof of Lemma~\ref{lem:smooth}. 

\section{Numerical results} \label{sec:results}

\begin{figure}[ht!]
 \centering
 \includegraphics[width=0.38\linewidth]{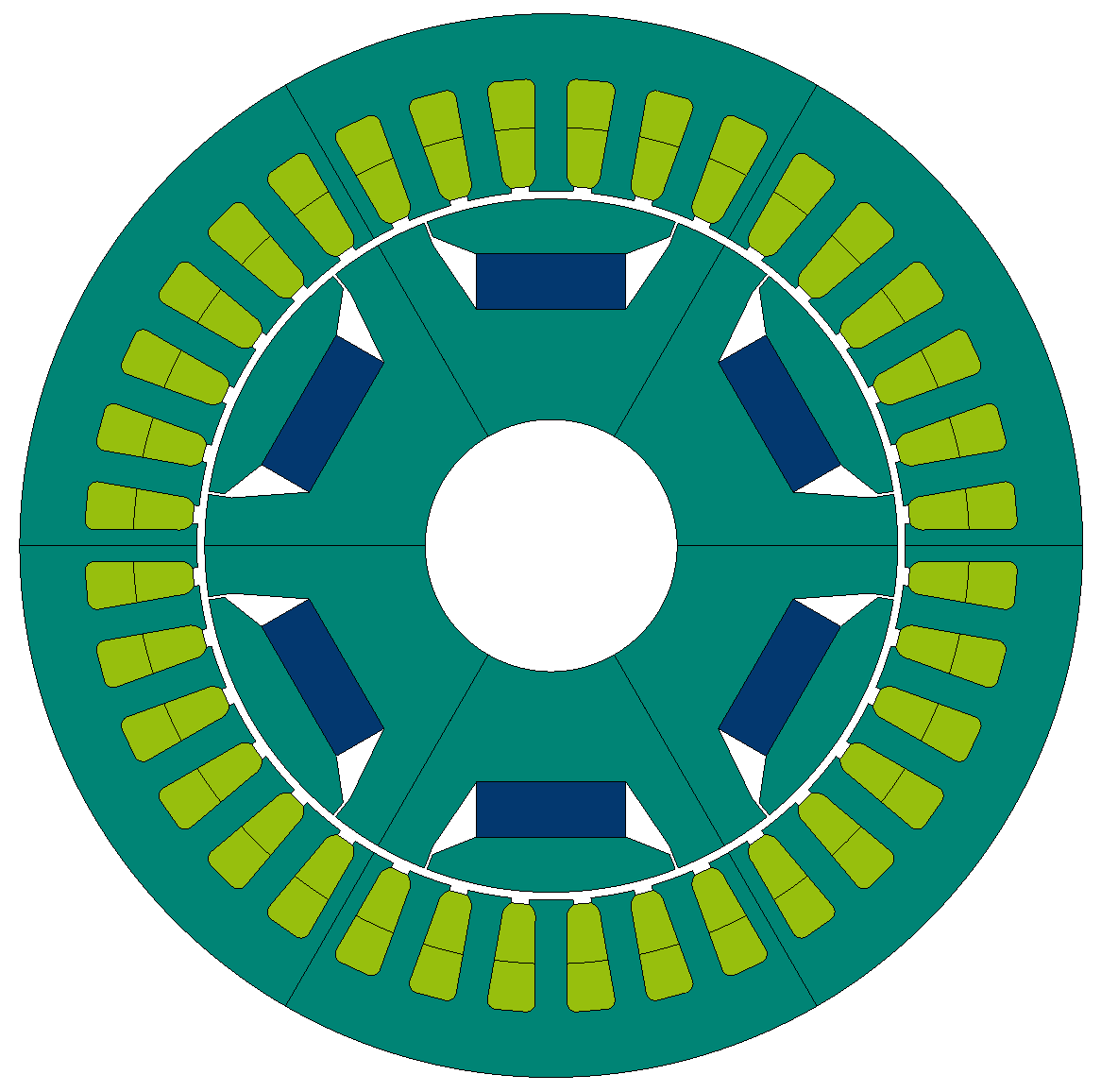}
 \includegraphics[width=0.61\linewidth]{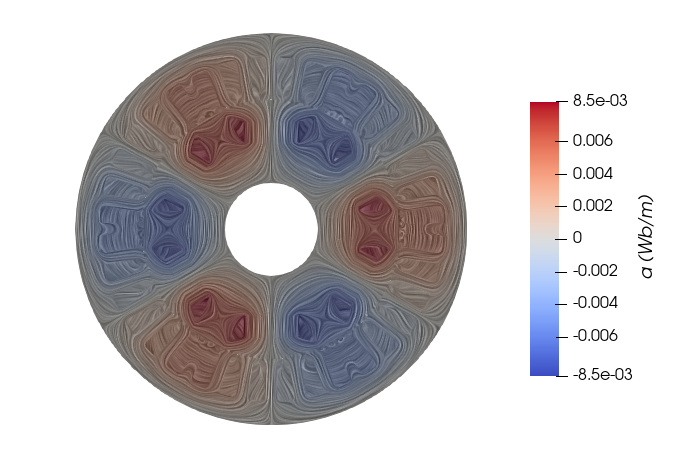}
 \caption{\textit{Left:} Model of a six-pole permanent magnet synchronous machine. The rotor on the inside consists of iron (green) and permanent magnets (blue). The stator on the outside consists of iron (green) and the slots for the copper windings (light green). \textit{Right:} Magnetic vector potential $a$ and magnetic field lines in the machine.}
 \label{fig:pmsm}
\end{figure}

In the following, we demonstrate the application of the harmonic mortar finite element and isogeometric analysis method for computing the torque in a six-pole permanent magnet synchronous machine shown in Fig.~\ref{fig:pmsm}. 
The geometry consists of a rotor with three pole pairs with embedded permanent magnets and a stator with $36$ slots filled with copper windings.
The length of the actual machine is $L=\SI{0.1}{m}$, so the energy and torque have to be scaled by $L$ as outlined above. 
Some information about the geometry and material parameters is summarized in Table~\ref{table:machine_parameters};
for details see \cite[Chapter V.A]{Bontinck_2018af}.
\begin{table}[ht!]
\centering
\caption{Geometry and material parameters of the machine.}
\label{table:machine_parameters}
\begin{tabular}{ll}
\hline
\multicolumn{2}{c}{geometry}\\
\hline
inner radius rotor & \SI{16}{mm}\\
outer radius rotor & \SI{44}{mm}\\
inner radius stator & \SI{45}{mm}\\
outer radius stator & \SI{67.5}{mm}\\
\hline
\multicolumn{2}{c}{materials}\\
\hline
relative permeability iron & 500\\
relative permeability copper & 1\\
relative permeability permanent magnets & 1.05\\
remanence field permanent magnet & \SI{0.94}{T}\\
\end{tabular}
\end{table}

\medskip 

\textbf{Simulation setup.}
In our numerical tests, we compute the cogging torque via formula \eqref{eq:torqueh2} in the absence of excitation currents $j_e=0$.
For the discretization of the magnetostatic problems in the stator and rotor domains, we used open source package \textsc{GeoPDEs} \cite{Vazquez_2016aa}. The computations are carried out with B-splines of degree $p=1$ and $p=2$ with maximal continuity as basis functions for the solution, which corresponds to a standard finite element method (FEM) and a multipatch spline discretization with patchwise $C^1$-continuous splines of second order (IGA). To accommodate for material discontinuities, only $C^0$-continuity is enforced across patch borders; see \cite{Bontinck18} for details. 
At the external boundaries of the stator and rotor, the magnetic vector potential $a$ is set to zero. 
For the computations of different degree $p$ the mesh refinement has been adapted such that the number of degrees of freedom is comparable for $p=1$ and $p=2$. The number of degrees of freedom for $p=1$ is $N_{\mathrm{dof},R}=26532$ in the rotor domain and $N_{\mathrm{dof},S}=39528$ in the stator domain. For $p=2$ the number of degrees of freedom is $N_{\mathrm{dof},R}=25140$ in the rotor domain and $N_{\mathrm{dof},S}=37764$ in the stator domain.
The reference solution is computed on a  refined mesh with $p=2$ leading to $N_{\mathrm{dof},R}=93360$ and $N_{\mathrm{dof},S}=128664$ degrees of freedom in the rotor and stator domain, and with trigonometric degree $N=200$. 
The saddlepoint problems \eqref{eq:lin1}--\eqref{eq:lin3} are solved efficiently by a Schur complement method;  details are given in Section~\ref{sec:appendix}.

\medskip

\textbf{Symmetry of the solution.}
A snapshot of the magnetic vector potential and its isolines, which correspond to the magnetic flux lines, is depicted in Fig.~\ref{fig:pmsm}.
Note that by symmetry of the rotor and stator geometry, the solution components $(a_S,a_R,\lambda)$ have a $60^\circ$ anti-symmetry. In particular, for the Fourier expansion of the computed Lagrange multiplier
\begin{align} \label{eq:lagrange}
    \lambda_N(\theta) = \frac{c_0}{2} + \sum_{n=1}^N c_n \cos (n \theta) + d_n \sin(n \theta),
\end{align}
only the coefficients $c_n$, $d_n$ with index $n\in \mathcal{N}=\{3,9,15,...\}$
are non-zero. 
By careful construction of the meshes, the geometric symmetries are preserved on the discrete level and the same behavior is observed for the numerical solution, as can be seen from Table~\ref{tab:lagrange}.  
\begin{table}[ht!]
\centering
\caption{Relevant and irrelevant Fourier modes of the Lagrange multipliers of an IGA simulation result for $\alpha=\SI{7}{\degree}$, $C_n=\sqrt{c_n^2+d_n^2}$ is the total amplitude.\label{tab:lagrange}}
\begin{tabular}{c|c|c|c}
$C_3$ & $C_9$ & $C_{15}$ & $\sum\nolimits_{n\not\in\mathcal{N}}C_n$ \\
\hline 
$\SI{9.2621e+03}{}$ & $\SI{1.6651e+03}{}$ &  $\SI{5.5848e+03}{}$ & $\SI{4.9720e-06}{}$ 
\end{tabular}
\end{table}
Using this prior knowledge, most of the Lagrange multipliers can be eliminated already during assembling leading to more efficient computations, in particular in connection with the Schur complement technique outlined in Section~\ref{sec:appendix}.

\medskip 

\begin{figure}[ht!]
\centering
\includegraphics{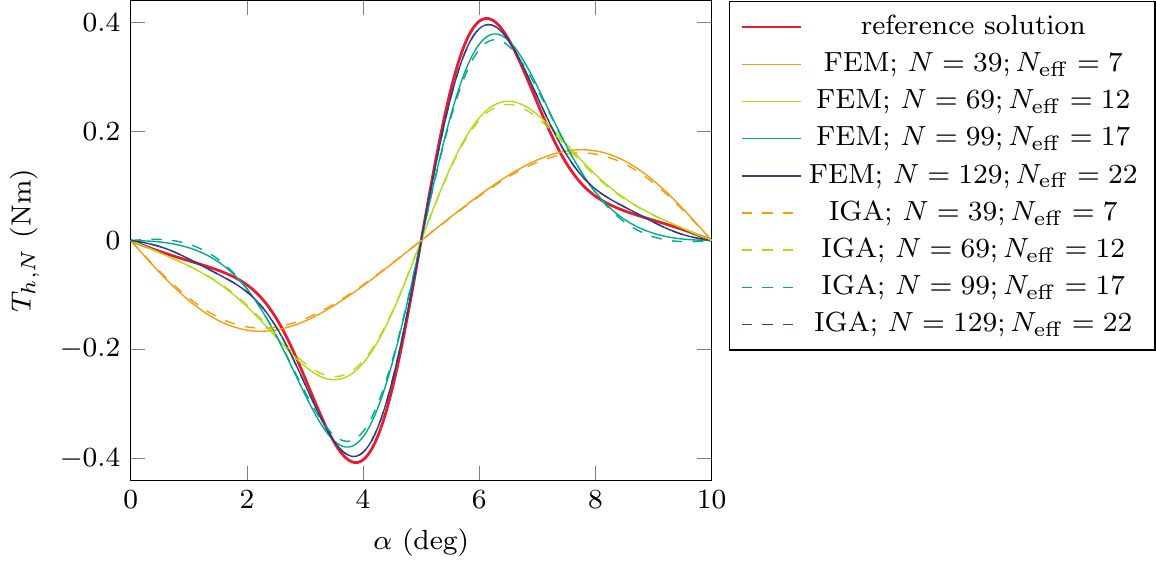}
\caption{Cogging torque of the machine. Computed using trigonometric polynomials of degree $N$ as basis for the Lagrange multiplier space. Computed with FEM and IGA using basis functions of degree $p=1,2$.}
 \label{fig:torque}
 \end{figure}

 \textbf{Torque computation.}
In Fig.~\ref{fig:torque}, we display the torque computed via formula \eqref{eq:torque} with FEM ($p=1$) and IGA ($p=2$) for different orders $N$ of harmonic basis functions. 
By the $10^\circ$ symmetry of the stator slots and the corresponding symmetry of the mesh, the computed torque can be shown to exhibit a periodic behavior with a period length of $10^\circ$. 
Moreover, the torque can be seen to be anti-symmetric with respect to half of the period. We further observe that the torque converges towards that of the reference solution when increasing the number $N$ of harmonic basis functions. As explained by the theoretical results in \cite{Egger_2021aa}, the computation becomes unstable when increasing the dimension of the Lagrange multiplier spaces too much; here for $N=250$. 

Due to its periodicity, the torque $T_{h,N}$ can be approximated by a Fourier series  
\begin{align} \label{eq:torque_fourier} 
T_{h,N}^M(\alpha) = \frac{\hat{c}_0}{2} + \sum\limits_{m=1}^{M} \hat{c}_m \cos(m \alpha) + \hat{d}_m \sin(m \alpha)
\end{align}
with coefficients $\hat{c}_m, \hat{d}_m \in\RR$.
By the symmetry properties discussed above, one can see that only $\hat{d}_m$ for $m\in \mathcal{M}=\{36,72,108,...\}$ is non-zero, while $\hat{d}_m=0$ of all $m \not\in \mathcal{M}$ and $\hat{c}_m=0$ for all $m$.
The corresponding results of our computations are shown in Table~\ref{table:torque_modes}.
\begin{table}[ht!]
\centering
\caption{Relevant and irrelevant Fourier modes of the torque of an IGA simulation result for a full rotation of the machine.}
\label{table:torque_modes}
\begin{tabular}{c|c|c|c|c}
$|\hat{d}_{36}|$ & $|\hat{d}_{72}|$ & $|\hat{d}_{108}|$ & $\sum\limits_{m}|\hat{c}_m|$ & $\sum\limits_{m\not\in\mathcal{M}}|\hat{d}_m|$ \\
\hline 
$\SI{2.293e-01}{}$ & $\SI{1.784e-01}{}$ &  $\SI{9.17e-02}{}$ & $\SI{6.2079e-11}{}$ & $\SI{51.7843e-11}{}$ 
\end{tabular}
\end{table}

\section{Discussion} \label{sec:conclusions}

We have discussed the harmonic mortar method using trigonometric polynomials as Lagrange multipliers for the efficient treatment of rotation in electric machine computation. The energy balance has been used to derive an explicit formulation for the computation of the torque in the continuous and discretized settings. It has been shown that the electro-magneto-mechanic energy balance holds exactly also on the discrete level.
Furthermore, if the functions in the Lagrange multiplier space are smooth then this is inherited by the torque.
The derived formulation has been applied to a typical permanent magnet synchronous machine model using isogeometric analysis and a lowest order finite element method with an exact geometry mapping. The cogging torque has been computed with the proposed method and the influence of the number of Lagrange multiplier basis functions on the solution has been investigated.
\section*{Acknowledgements}
This work is supported by the Graduate School CE within the Centre for Computational Engineering at Technische Universität Darmstadt, by the German Research Foundation via the project SCHO 1562/6-1, and by the Defense Advanced Research Projects Agency (DARPA), under contract HR0011-17-2-0028. The views, opinions and/or findings expressed are those of the author and should not be interpreted as representing the official views or policies of the Department of Defense or the U.S. Government.

\appendix 

\section{Appendix}
\label{sec:appendix}
If the machine model uses only a small number of interface modes $N$ then local condensing is very effective. Let as assume that the matrices $K_S$ and $K_R$ are non-singular, e.g. due to Dirichlet boundaries of the stator and rotor domains, then the system \eqref{eq:lin1}-\eqref{eq:lin3} can be rewritten as
 \begin{alignat}{5}
    a_S(\alpha) & &    &+K_S^{-1}B_S^\top \lambda(\alpha) &&= K_S^{-1}j_e \\
    a_R(\alpha) & &    &-K_R^{-1}B_R(0)^\top R(\alpha)^\top \lambda(\alpha) \ &&= K_R^{-1}j_M \\
    B_S a_S(\alpha) &-& R(\alpha)B_R(0) a_R(\alpha)& &&= 0. 
\end{alignat}
The internal degrees of freedom are straightforwardly eliminated using the Schur-complement. This gives rise to the low-dimensional interface problem
\begin{align}
{K}_{\textrm{int}}(\alpha)
{\lambda}
&=
{f}_{\textrm{int}}(\alpha) \label{eq:sysint}
\end{align}
with
\begin{align}
{K}_{\textrm{int}}(\alpha)&=B_SK_S^{-1}B_S^\top+R(\alpha)B_R(0)K_R^{-1}B_R(0)^\top R(\alpha)^\top, \label{eq:kint}
	\\
	{f}_{\textrm{int}}(\alpha)&=B_SK_S^{-1}j_e-R(\alpha)B_R(0)K_R^{-1}j_M. \label{eq:fint} 
\end{align}
The inverses in \eqref{eq:kint} and \eqref{eq:fint} are not needed explicitly. Instead, one factorization and a few forward/backward substitutions (for each spectral basis at the interface) can be used to precompute the necessary expressions.
Thus, only the small system \eqref{eq:sysint} has to be solved in the online phase for different rotation angles. 

The internal degrees of freedom can be cheaply reconstructed, i.e.,
\begin{align}
	a_S(\alpha)&=K_S^{-1}j_e-K_S^{-1}B_S^\top \lambda(\alpha),\label{eq:recostructionust}\\
	a_R(\alpha)&=K_R^{-1}j_M+K_R^{-1}B_R(0)^\top R(\alpha)^\top \lambda(\alpha).\label{eq:recostructionurt}	
\end{align}
This is computationally convenient when dealing with rotation since only the low-dimensional matrix ${R}(\alpha)$ depends on $\alpha$. This leads to a significant reduction in, e.g., when considering the computation of a rotating electric machine.


\begin{thebibliography}{10}

\bibitem{Alotto01}
P.~Alotto{\em ~et~al}.
\newblock Discontinuous finite element methods for the simulation of rotating
  electrical machines.
\newblock {\em COMPEL}, 20:448--462, 2001.

\bibitem{Buffa01}
A.~Buffa, Y.~Maday, and F.~Rapetti.
\newblock A sliding mesh-mortar method for a two dimensional currents model of
  electric engines.
\newblock {\em ESAIM Math. Model. Numer. Anal.}, 35:191--228, 2001.

\bibitem{Lange10}
E.~Lange, F.~Henrotte, and K.~Hameyer.
\newblock A variational formulation for nonconforming sliding interfaces in
  finite element analysis of electric machines.
\newblock {\em IEEE Trans. Magn.}, 46:2755--2758, 2010.

\bibitem{Tsukerman95}
I.~Tsukerman.
\newblock Accurate computation of 'ripple solutions' on moving finite element
  meshes.
\newblock {\em IEEE Trans. Magn.}, 31:1472--1475, 1995.

\bibitem{HenrotteHameyer04}
F.~Henrotte and K.~Hameyer.
\newblock Computation of electromagnetic force densities: Maxwell stress tensor
  vs. virtual work principle.
\newblock {\em J. Comput. Appl. Math.}, 168:235--243, 2004.

\bibitem{Bontinck18}
Z.~Bontinck{\em ~et~al}.
\newblock Isogeometric analysis and harmonic stator-rotor coupling for
  simulating electric machines.
\newblock {\em Comput. Meth. Appl. Mech. Engrg.}, 334:40--55, 2018.

\bibitem{DeGersem04}
H.~De~Gersem and T.~Weiland.
\newblock Harmonic weighting functions at the sliding interface of a
  finite-element machine model incorporating angular displacement.
\newblock {\em IEEE Trans. Magn.}, 40:545--548, 2004.

\bibitem{Bossavit90}
A.~Bossavit.
\newblock Forces in magnetostatics and their computation.
\newblock {\em J. Appl. Phys.}, 67:5812--5814, 1990.

\bibitem{Coulomb83}
J.~L. Coulomb.
\newblock A methodology for the determination of global electromechanical
  quantities from a finite element analysis and its application to the
  evaluation of magnetic forces, torques and stiffness.
\newblock {\em IEEE Trans. Magn.}, 16:2514--2519, 1983.

\bibitem{Brezzi74}
F.~Brezzi.
\newblock On the existence, uniqueness and approximation of saddle-point
  problems arising from {L}agrangian multipliers.
\newblock {\em RAIRO Anal. Numer.}, 8:129--151, 1974.

\bibitem{Egger_2021aa}
H.~Egger{\em ~et~al}.
\newblock On the stability of harmonic mortar methods with application to
  electric machines.
\newblock In W.~Schilders, M.~van Beurden, and N.~Budko, editors, {\em
  Scientific Computing in Electrical Engineering {SCEE} 2020}, Mathematics in
  Industry, Berlin, Springer, 2021.

\bibitem{Braess}
D.~Braess.
\newblock {\em Finite Elements}.
\newblock Cambridge University Press, New York, 3rd edition, 2007.

\bibitem{Cottrell_2009aa}
J.~A. Cottrell, T.~J.~R. Hughes, and Y.~Bazilevs.
\newblock {\em Isogeometric Analysis: Toward Integration of {CAD} and {FEA}}.
\newblock Wiley, 2009.

\bibitem{Neuman_1981aa}
E.~Neuman.
\newblock Moments and fourier transforms of b-splines.
\newblock {\em J. Comput. Appl. Math.}, 7(1):51--62, 1981.

\bibitem{Bontinck_2018af}
Z.~Bontinck.
\newblock {\em Simulation and Robust Optimization for Electric Devices with
  Uncertainties}.
\newblock Dissertation, Technische Universität Darmstadt, 2018.
\newblock \url{https://tuprints.ulb.tu-darmstadt.de/id/eprint/8330}.

\bibitem{Vazquez_2016aa}
R.~Vázquez.
\newblock A new design for the implementation of isogeometric analysis in
  {Octave} and {Matlab}: {GeoPDEs} 3.0.
\newblock {\em Comput. Math. Appl.}, 72:523--554, 2016.
  \href{https://dx.doi.org/10.1016/j.camwa.2016.05.010}{doi:10.1016/j.camwa.2016.05.010}.

\end{thebibliography}
\end{document}